\documentclass{amsart}

\usepackage{amssymb}
\usepackage{lipsum}

\usepackage{ytableau}
\usepackage{amsmath}

\usepackage[hyphens]{url}
\usepackage[colorlinks=true, allcolors=blue, breaklinks=true]{hyperref}

\begin{document}

		\begin{center}
			{\Large\bf Universal quantum dimensions II: 
				$\gamma$-dependent factors }\\
			\vspace*{1 cm}
			{\large  }
			\vspace*{0.5 cm}
			
			{\large  M.Y. Avetisyan and R.L.Mkrtchyan 
			}

			\vspace*{0.5 cm}
			
			{\small\it Alikhanyan National Science Laboratory (Yerevan Physics Institute), \\ 2 Alikhanian Br. Str., 0036 Yerevan, Armenia}

			\vspace*{0.5 cm}

		\end{center}

		{\bf Abstract.}
		
			We develop an approach for computing the universal, in Vogel's sense, quantum dimensions. In a previous work a method for calculation of $\gamma$-independent part of these dimensions for universal multiplets with nonzero  associated members for classical algebras was proposed. In the present paper, this approach is extended to the $\gamma$-dependent factors. In particular, we derive these contributions for the $E$ universal multiplet, which appears in the universal decomposition of the fourth power of the adjoint representation, thereby obtaining its universal quantum dimension for the first time.
			
			The analysis requires knowledge of the quantum dimensions of the classical and $E_8$ members of the multiplet; for the remaining exceptional algebras, the corresponding quantum dimensions are recovered automatically, providing an additional consistency check of the approach. We also partially extend a previously conjectured relation between the $sl$ and $so, sp$ members - formulated in terms of vertical and horizontal sums of Young diagrams - to the exceptional case. Finally, we present an explicit algorithm implementing the proposed method, which enables the derivation of universal formulae for higher-dimensional representations. The  universal quantum dimensions provide, in some cases, a unified description of  Chern–Simons observables, including knot invariants and partition functions, in a form valid simultaneously for all gauge algebras.

		

	
	

	\section{Introduction}

	Since Vogel's paper \cite{Vogel1999}, see also \cite{Vogel2011}, the "universality" approach in the theory and applications of simple Lie algebras has been actively developed. This approach reveals hidden connections between algebras by expressing a number of quantities (dimensions of irreducible representations, eigenvalues of Casimir operators, partition functions, etc.) in a universal form, namely in terms of universal (Vogel's) parameters $\alpha, \beta, \gamma$. These parameters are in fact homogeneous coordinates on the projective plane, and by specializing them to the points from Vogel's table \ref{tab:Vogel} one obtains the corresponding result for a particular simple Lie algebra. Determining the complete list of such universal quantities, as well as systematic methods for deriving them, remains an important open problem.
	
	As an example of universal representation show the free energy of Chern-Simons theory on $S^3$. It is calculated by \cite{Witten1989} and transformed into the following universal form in \cite{KreflMkrtchyan2015}:

	\begin{eqnarray}\label{Free}
		F=\frac{1}{4}\int_{R_+} \frac{dx}{x} \frac{\sinh\left(x(t-\delta)\right)}{\sinh\left( x t \right)\sinh\left(x \delta\right)} f(2x) \\ \nonumber
			\delta=\kappa+t, \,\,\,t=\alpha+\beta+\gamma
	\end{eqnarray}
	where $\kappa$ is the coupling of Chern-Simons theory, and $f(x)$ is the universal form of quantum dimension of adjoint representations, derived in \cite{Westbury2003}: 
	\begin{align}\label{cad}
		f(x)&=\frac{\sinh(x\frac{\alpha-2t}{4})}{\sinh(x\frac{\alpha}{4})}\frac{\sinh(x\frac{\beta-2t}{4})}{\sinh(x\frac{\beta}{4})}\frac{\sinh(x\frac{\gamma-2t}{4})}{\sinh(x\frac{\gamma}{4})} 	
	\end{align}
	Substituting into \eqref{Free} values of universal parameters from table \ref{tab:Vogel} we obtain the Chern-Simons free energy for the corresponding gauge algebra. 
	
	In this paper we develop an approach for calculating universal quantum dimensions of a certain class of universal multiplets, namely those with nonzero {\it associated} members for classical algebras. In \cite{Mkrtchyan2026} the first part of this method was developed, for the calculation of factors independent of $\gamma$ (which is the only universal parameter for classical algebras in the table \ref{tab:Vogel} that depends on the rank). There it was used to determine such terms for the universal multiplet $E$, considered in the \cite{CohenMan1996,AvetisyanIsaevKrivonosMkrtchyan2024}. In the present paper we consider   $\gamma$-dependent factors of universal multiplets and particularly obtain the remaining $\gamma$-dependent factors  and hence the final universal quantum dimension formula for the multiplet $E$.
	
	There exist several approaches for deriving universal formulae for dimensions or quantum dimensions of representations: \cite{Vogel1999,Vogel2011,LandsbergManivel2006,Mkrtchyan2017,IsaevKrivonosProvorov2023}. However, none of them is presently capable of handling arbitrary representations. The main advantage of the method proposed in \cite{Mkrtchyan2026} and further developed in the present paper is that it is expected to be applicable to arbitrary "eligible" universal multiplets, namely those with nonzero associated members for both classical and exceptional algebras.

	We call the {\it universal multiplet} corresponding to a given universal (quantum) dimension formula (sometimes referred to as a {\it big universal multiplet}) the union of {\it small universal multiplets}. The latter are sets of representations of a given simple Lie algebra obtained by evaluating the universal dimension formula at the values of Vogel's parameters corresponding to that algebra in table \ref{tab:Vogel}, together with their permutations. We call {\it associated} the set of representations of different algebras (one from each algebra), obtained from the same formula by substituting values of Vogel's parameters from table \ref{tab:Vogel}. Similarly, associated are representations, obtained from the given formula with  simultaneously permuted parameters from table  \ref{tab:Vogel}.
	
	The idea of the method is based on the observation of a correspondence between associated representations of classical algebras \cite{Mkrtchyan2025b,Mkrtchyan2026}, which is partially extended here to exceptional algebras. It is assumed (and verified in several examples) that this correspondence induces relations between the arguments of the hyperbolic sines appearing in the quantum dimensions, allowing one to reconstruct their universal form uniquely. This situation is similar to the phenomenon whereby universal eigenvalues of the second-order Casimir operator, or of higher Casimir operators, can be uniquely reconstructed from their values on classical algebras. In this context it is important to take into account certain subtleties related to automorphism groups of Dynkin diagrams, in particular the $Z_2$ automorphism of the $sl$ algebra. Throughout the paper, by the Casimir operator we always mean the quadratic Casimir operator.
	
	As an example, consider the simplest universal quantum dimension formula, namely that for the adjoint representation \eqref{cad}. 

	Its small universal multiplets consist of a single representation, namely the adjoint representation for the given algebra, since the formula is invariant under permutations. To illustrate the approach, consider the $\gamma$-independent terms in (\ref{cad}) and attempt to derive them from the quantum dimensions of the adjoint representations of classical algebras given below.
	
	For $sl(N)$:
	\begin{eqnarray}\label{qsl}
		\frac{\sinh\left(\frac{x}{2}(N-1)\right) \sinh\left(\frac{x}{2}(N+1)\right)}{\sinh\left(\frac{x}{2}\right)\sinh\left(\frac{x}{2}\right)}
	\end{eqnarray}
	
	For $so(N)$:
	
	\begin{eqnarray}\label{qso}
		\frac{\sinh\left(\frac{x}{4}N\right) \sinh\left(\frac{x}{2}(N-1)\right) \sinh\left(\frac{x}{2}(N-4)\right)}{\sinh\left(\frac{x}{2}\right)\sinh\left(x\right)\sinh\left(\frac{x}{4}(N-4)\right) }
	\end{eqnarray}
	
	For $sp(N)$:
	
	\begin{eqnarray}\label{qsp}
		\frac{\sinh\left(\frac{x}{8}N\right) \sinh\left(\frac{x}{4}(N+1)\right) \sinh\left(\frac{x}{4}(N+4)\right)}{\sinh\left(\frac{x}{4}\right)\sinh\left(\frac{x}{2}\right)\sinh\left(\frac{x}{8}(N+4)\right) }
	\end{eqnarray}
	
	As can be seen from these three expressions, the $\gamma$-independent part of the universal expression must contain at least two hyperbolic sines in the denominator. We parameterize them as follows:
	
	\begin{align}\label{uqs}
		&\sinh\left(\frac{x}{4}(x_1 w_x+y_1 w_y)\right)\sinh\left(\frac{x}{4}(x_2 w_x+y_2 w_y)\right)  \\
		&   w_x=-\frac{2\alpha+\beta}{2}, \quad		w_y=\frac{\alpha+\beta}{2}
	\end{align}
	where we introduced a convenient parameterization in terms of $w_x$ and $w_y$, which for classical algebras take the values listed in table \ref{tab:wxwy}.

	\begin{table}[ht] 
		\centering    
		\caption{Values of $w_x, w_y$}
		\label{tab:wxwy}
		\begin{tabular}{|c|c|c|c|}
			\hline
			& $s$l & $so$ & $sp$ \\
			\hline
			$w_x$ & 1 & 0 & $\frac{3}{2}$ \\
			\hline
			$w_y$ & 0 & 1 & $-\frac{1}{2}$ \\
			\hline
		\end{tabular}
	\end{table}

	Now compare (\ref{uqs}) with the denominator of the quantum dimension for $sl$ (\ref{qsl}), assuming all arguments are positive. It follows that $x_1=x_2=2$. Comparison with the $so$ formula (\ref{qso}) leads to either $y_1=2, y_2=4$ or vice versa, which is inessential (since the two possibilities are related by renumeration), because $x_1=x_2$. From these comparisons we deduce the final form of the denominator
	
	\begin{align}\label{uqad} \nonumber
	&	\sinh\left(\frac{x}{4}(x_1 w_x+y_1 w_y)\right)\sinh\left(\frac{x}{4}(x_2 w_x+y_2 w_y)\right) = \\
	&	\sinh\left(\frac{x}{4}(-\alpha)\right)\sinh\left(\frac{x}{4}(\beta)\right)
	\end{align}
	in agreement with (\ref{cad}). Note that the $sp$ data were not used in deriving the complete universal answer; however, they are naturally reproduced by the final result, as well as the cases of exceptional algebras.
	
	Next we turn to the $\gamma$-dependent terms, which constitute the main subject of the present paper. From the expressions above for classical algebras we assume that such terms consist of $2\gamma$-terms, i.e. those of the type $\sinh\left(\frac{x}{4}(2\gamma+...)\right)$, of which there should be at least two in the numerator, and $1\gamma$-terms of the type $\sinh\left(\frac{x}{4}(\gamma+...)\right)$, which must appear in equal numbers (at least one) in the numerator and denominator.
	
	Assuming the minimal number of $2\gamma$-terms — two in the numerator — we parameterize them as
	
	\begin{align}\label{2gad}
		\sinh\left(\frac{x}{4}(2\gamma+x_1 w_x+y_1 w_y)\right)\sinh\left(\frac{x}{4}(2\gamma+x_2 w_x+y_2 w_y)\right)
	\end{align}
	
	One can easily verify, by comparison with (\ref{qsl}) and (\ref{qso}), that $\{x_1,x_2\}=\{-2,2\}$ and $\{y_1,y_2\}=\{0,6\}$, where unordered sets are understood. One may fix the numbering by choosing $x_1=-2, x_2=2$, while the numbering of $y_i$ remains undetermined. Of course, one could check both possibilities, but such an approach becomes impractical for higher orders. The correct assignment of $y_i$ can be established by comparison with the $sp$ result (\ref{qsp}). In this case (\ref{2gad}) becomes
	
	\begin{align}\label{2gadsp}
		\sinh\left(\frac{x}{4}(N+4+x_1 \frac{3}{2}-y_1 \frac{1}{2})\right)\sinh\left(\frac{x}{4}(N+4+x_2 \frac{3}{2}-y_2 \frac{1}{2})\right)
	\end{align}
	
	Comparing with (\ref{qsp}), we obtain $\{1-\frac{1}{2}y_1, 7-\frac{1}{2}y_2\}=\{1,4\}$, from which (recalling $\{y_1,y_2\}=\{0,6\}$) the unique solution $y_1=0, y_2=6$ follows. This gives the correct form of the $2\gamma$-terms and provides an additional consistency check. The $1\gamma$-terms can be recovered in a similar way.
	
	In the following sections we perform a more involved analysis for the quantum dimension of the universal multiplet $E$. In the previous work the $\gamma$-independent terms were derived; here we determine the $\gamma$-dependent terms and obtain the final universal expression reproducing the quantum dimensions of all representations in the universal multiplet, as well as perform some checks. In the Conclusion we present the method in its general form.

	\section{$E$ universal multiplet}
	
	The universal dimension of the $E$ universal multiplet was given in \cite{AvetisyanIsaevKrivonosMkrtchyan2024}:

		\begin{align}\label{dimE}  \nonumber
		\dim E&=
		-\frac{64 (\alpha+\gamma) (2 \alpha+\gamma) (\alpha+2 \gamma) (\beta+\gamma) (2 \beta+\gamma) }{\alpha^2 \beta^2 \gamma (\alpha-3 \beta) (\alpha-\beta)^2 (3 \alpha-\beta) (\alpha-\gamma) } \times \\ \nonumber 
		&\frac{  (\beta+2 \gamma) (\alpha+\beta+\gamma)	(2 \alpha+\beta+\gamma) 	(\alpha+2 \beta+\gamma) }{(\beta-\gamma) (\alpha+\beta-2 \gamma)} \times \\
&(2 \alpha+2 \beta+\gamma) (2 \alpha+\beta+2 \gamma) (\alpha+2 \beta+2 \gamma)
	\end{align}

	The only nonzero member of the $sl$ small multiplet of $E$ is the one obtained by substituting into (\ref{dimE}) the parameters from table \ref{tab:Vogel}:
	
	\begin{align}\label{slE} \nonumber
		&D_s(\ydiagram{3,1},\ydiagram{2,1,1})\equiv \\
		&(210...0101)_s\equiv (210...0101)\oplus(101...012)
	\end{align}
	where we use the notations $D_s(?,?)$ and explanation of that notation in terms of Dynkin labels. Note that this is the sum of two irreducible representations of $sl$, which becomes irreducible upon extending the group by the $Z_2$ automorphism of its Dynkin diagram. 
	
	The dimension of this representation is twice the dimension of each irrep in (\ref{slE}) and is given by
	\begin{align}
		\frac{1}{32} (N-4) (N-2) (N-1)^2 (N+1)^2 (N+2) (N+4)
	\end{align}
	It is invariant under $N \leftrightarrow -N$ because of the symmetry of (\ref{dimE}) under $\alpha \leftrightarrow \beta$. The other two members of the $sl$ small multiplet vanish. 
	
	The associated $so$ representation is given by the vertical sum (\cite{Mkrtchyan2025b}) of the $sl$ diagrams in (\ref{slE}):
	\begin{eqnarray}\label{soE}
		\ydiagram{3,2,1,1,1}
	\end{eqnarray}
	Its dimension is
	\begin{align}
		\frac{1}{630} (N-6) (N-4) (N-2) (N-1) N^2 (N+2) (N+4) 
	\end{align}
	while the other elements of the $so$ small multiplet are zero. 
	
	The vertical sum of two Young diagrams $\lambda,\tau$ is another Young diagram $\lambda \oplus_v \tau$  defined as follows.  Let $h^\mu_i$ be the height of the $i$-th column of the Young diagram $\mu$, numerated from left to right.  	Then
	\begin{eqnarray}
		h^{\lambda \oplus_v \tau}_i = h^\lambda_i + h^\tau_i .
	\end{eqnarray}
	
	The associated $sp$ representation is given by the horizontal sum (\cite{Mkrtchyan2025b,Mkrtchyan2026}) of the $sl$ diagrams in (\ref{slE}):
	
	\begin{eqnarray}\label{spE}
		\ydiagram{5,2,1}
	\end{eqnarray}
	with dimension
	\begin{eqnarray} 
		\frac{1}{630} (N-4) (N-2) N^2 (N+1) (N+2) (N+4) (N+6) 
	\end{eqnarray}
	The other two members of the $sp$ small multiplet are zero. 
	
		The horizontal sum of two Young diagrams is defined similarly to vertical one. For two Young diagrams $\lambda,\tau$ the horizontal sum is another Young diagram $\lambda \oplus_h \tau$  defined as follows.  Let $r^\mu_i$ be the length of the $i$-th row of the Young diagram $\mu$, numerated from up to down.	Then
	\begin{eqnarray}
		r^{\lambda \oplus_h \tau}_i = r^\lambda_i + r^\tau_i .
	\end{eqnarray}
	
	The quantum dimensions of the representations in (\ref{slE}) are equal, due to the $Z_2$ symmetry, and are given by

	\begin{align}\label{qsle} \nonumber
	&	\frac{
			\sinh\!\big(\frac{x}{2}(N-4)\big)\,\sinh\!\big(\frac{x}{2}(N-2)\big)\,\sinh\!\big(\frac{x}{2}(N-1)\big)^{2}\,			
		}{
			\sinh(\frac{x}{2})^{4}\,\sinh(x)^{2}\,\sinh(2x)^{2} 
		} \times \\  
&	\sinh\!\big(\frac{x}{2}(N+1)\big)^{2}\,	\sinh\!\big(\frac{x}{2}(N+2)\big)\,\sinh\!\big(\frac{x}{2}(N+4)\big)  
	\end{align}

	The quantum dimension of the $so(2n)$ representation (\ref{soE}) is

		\begin{align}\label{qsoe} \nonumber
		&\frac{
			\sinh(\frac{x}{2}n)\,\sinh(\frac{x}{2}(n+2))\,
			\sinh(x(n-5))\,
		}{
			\sinh(\frac{x}{2})^{3}\,\sinh(x)\,\sinh(3\frac{x}{2})^{2}\,\sinh(5\frac{x}{2})\,
		} \times \\  \nonumber
		&\frac{\sinh(x(n-3))\,\sinh(x(n-2))\, 	\sinh(x(n-1))^{2}\, }{\sinh(7\frac{x}{2})\,\sinh(\frac{x}{2}(n-5))\,\sinh(\frac{x}{2}(n-1))} \times \\
	&	\sinh(\frac{x}{2}(2n-1))\,\sinh(x n)\,\sinh(x(n+1))
	\end{align}

	Similarly, the quantum dimension of the $sp(2n)$ representation (\ref{spE}) is

		\begin{align}\label{qspe} \nonumber
		&\frac{
			\sinh(\frac{x}{4}(n-2))\,\sinh(\frac{x}{4}n)\,\sinh(\frac{x}{2}(n-1))\,\,
		}{
			\sinh(\frac{x}{4})^{3}\,\sinh(\frac{x}{2})\,\sinh(3\frac{x}{4})^{2}\,\sinh(5\frac{x}{4})\,
		}  \times\\  \nonumber
		&\frac{  \sinh(\frac{x}{2}(n))\,\sinh(\frac{x}{4}(2n+1))  \sinh(\frac{x}{2}(n+1))^{2}\,}{  \sinh(7\frac{x}{4})\,\sinh(\frac{x}{4}(n+1))\,  \sinh(\frac{x}{4}(n+5))} \times \\
		&\sinh(\frac{x}{2}(n+2))\,\sinh(\frac{x}{2}(n+3))\,\sinh(\frac{x}{2}(n+5))
	\end{align}

	The $\gamma$-independent factor $Q_0$ of the universal quantum dimension of the representation $E$ is (\cite{Mkrtchyan2026}):

	\begin{align}\label{udenomE2}\nonumber
		Q_0 &=	\frac{\sinh(x(\alpha+\beta)/2)}{\sinh(x(\alpha+\beta)/4)}   \frac{1}{	\sinh\left(\frac{x}{4}\alpha\right)^2 	\sinh\left(\frac{x}{4}\beta\right)^2 }   \times \\ 
		&\frac{1}{	\sinh\left(\frac{x}{4}(\beta-\alpha)\right)^2 	\sinh\left(\frac{x}{4}(\beta-3\alpha)\right)\sinh\left(\frac{x}{4}(3\beta-\alpha)\right)}
	\end{align}

	We now turn to the $\gamma$-dependent terms. Recall that they consist of terms of the following type (for the $sl(N)$ and $so(N)$ algebras):
	\begin{eqnarray}
		\sinh\big(\frac{x}{4}(2N+const)\big)
	\end{eqnarray}	
	which we call $2\gamma$-terms, since their universal form is $\sinh(\frac{x}{4}(2\gamma+...))$, and 
	\begin{eqnarray}
		\sinh\big(\frac{x}{4}(N+const)\big)
	\end{eqnarray}	
	which we call $1\gamma$-terms. 
	
	Similar terms appear for the $sp(N)=sp(2n)$ algebra, with the rescaling $x\rightarrow x/2$; that is, the $2\gamma$-terms are $	\sinh\big(\frac{x}{2}(n+const)\big)$ and the $1\gamma$-terms are $	\sinh\big(\frac{x}{4}(n+const)\big)$.
	
	Let us now consider the $2\gamma$-terms in the expected universal quantum dimension formula for the universal multiplet E. In the quantum dimensions for $sl$, $so$, and $sp$, there are eight $2\gamma$-terms in the numerator and none in the denominator, so we first attempt to recover them using the same pattern of hyperbolic sines with universal arguments.
	
	Thus, we assume that the universal form of the $2\gamma$-terms is
	\begin{eqnarray}\label{uq8}
		\prod_{i=1}^{8} \sinh \left(\frac{x}{4}(2\gamma+x_i w_x+y_i w_y)\right)
	\end{eqnarray}
	
	Specializing this expression to $sl(N)$ and comparing with the quantum dimension (\ref{qsle}), we obtain the following values of $x_i$ (with the numbering fixed below):
	\begin{align}\label{xq8} \nonumber
		&x_1=-8,x_2=8,x_3=-4,x_4=4, \\
		&x_5=-2,x_6=2,x_7=-2,x_8=2
	\end{align}
	
	Similarly, specializing (\ref{uq8}) to $so(N)$ and comparing with (\ref{qsoe}), we obtain:
	\begin{eqnarray}\label{yq8}
		\{y_i, i=1...8\}=\{-12,-4,0,4,4,6,8,12\}
	\end{eqnarray}
	We stress that at this stage we do not know the individual values of the $y_i$; only the set of their values is determined. 
	
	Finally, specializing (\ref{uq8}) to $sp(N)$ and comparing with (\ref{qspe}), we obtain:
	\begin{eqnarray}\label{xyq8}
		\{\frac{3}{2}x_i-\frac{1}{2}y_i, i=1...8\}=\{-6,-4,-3,-2,-2,0,2,6\}
	\end{eqnarray}
	Again, only the set of values is known. Our task is therefore to determine specific values of $y_i$ satisfying (\ref{xq8}), (\ref{yq8}), and (\ref{xyq8}). 
	
	Consider first of these quantities, for  $i=1$, it is: 
	\begin{eqnarray}\label{11} 
		\frac{3}{2}x_1-\frac{1}{2}y_1=-12-\frac{1}{2}y_1. 
	\end{eqnarray}

	From the list (\ref{yq8}) of possible values of $y_1$ we see that the maximal value of   (\ref{11}) is $-6$, attained at $y_1=-12$. Since this is the minimal possible value from the list (\ref{xyq8}), we conclude that  $y_1=-12$.

	Consider the case  $i=2$ in (\ref{xyq8}): 
	\begin{eqnarray}
		\frac{3}{2}x_2-\frac{1}{2}y_2=12-\frac{1}{2}y_2
	\end{eqnarray}
	
	The minimal value of this quantity is 6, attained at $y_2=12$.  From the other hand, the maximal value in (\ref{xyq8}) is 6, so we obtain  $y_2=12$.

	Proceeding similarly, one obtains the unique solutions 
	\begin{eqnarray}
		y_3&=&-4, \quad \frac{3}{2}x_3-\frac{1}{2}y_3=-4, \\
		y_4&=&8, \quad \frac{3}{2}x_4-\frac{1}{2}y_4=2, \\
		y_5&=&0, \quad \frac{3}{2}x_5-\frac{1}{2}y_5=-3, \\
		y_6&=&6, \quad \frac{3}{2}x_6-\frac{1}{2}y_6=0
	\end{eqnarray}
	
	However, for $x_7=-2$ and $x_8=2$, the remaining sets $\{y_7,y_8\}=\{4,4\}$ and $\{\frac{3}{2}x_i-\frac{1}{2}y_i, i=7,8\}=\{-2,-2\}$ are incompatible. We therefore pass to the next possible minimal number of $2\gamma$-terms, namely nine hyperbolic sines in the numerator and one in the denominator, with the latter cancelling after specialization to one of the classical algebras. As we shall see, the previous calculations remain useful in this more general setting. 
	
	We now assume the following form for the $2\gamma$-terms in the universal quantum dimension formula for the universal multiplet $E$:
	\begin{eqnarray}\label{uq9} \nonumber
		\frac{\sinh(\frac{ x}{4}(2\gamma+\bar{x}w_x+\bar{y}w_y))}{\sinh(\frac{x}{4}(2\gamma+\bar{x}'w_x+\bar{y}'w_y))}	\times \\
		\prod_{i=1}^{8} \sinh \left(\frac{x}{4}(2\gamma+x_i w_x+y_i w_y)\right)
	\end{eqnarray}

	When specialized to $sl$, the hyperbolic sine in the denominator of this expression must cancel, say, the hyperbolic sine in the numerator of the first fraction. This implies that $\bar{x}'=\bar{x}$. The remaining factors must then reproduce the quantum dimension for $sl$. Thus, as before, the set of values of $x_i$ is again given by (\ref{xq8}), although at this point we do not yet fix their individual ordering:
	\begin{eqnarray}\label{xq9}
		\{x_i, i=1,...,8\}=\{-8,8,-4,4,-2,2,-2,2\}
	\end{eqnarray}
	
	Next, we specialize (\ref{uq9}) to $so$. The $2\gamma$-term in the denominator must now cancel one of the factors with $i=1,...,8$ in the numerator. (If it canceled the same numerator factor  in the first fraction, then that fraction would be identically equal to 1.) Suppose it cancels the eighth factor. This means that $\bar{y}'=y_8$. The remaining sine arguments must then coincide with those obtained earlier from (\ref{qsoe}):
	\begin{eqnarray}\label{yq9}
		\{y_i, i=1...7, \bar{y}\}=\{-12,-4,0,4,4,6,8,12\}
	\end{eqnarray}
	
	The $2\gamma$-term (\ref{uq9}) then takes the form
	\begin{align}\label{uspq9} \nonumber
		&\frac{\sinh(\frac{ x}{4}(2\gamma+\bar{x}w_x+\bar{y}w_y))}{\sinh(\frac{x}{4}(2\gamma+\bar{x}w_x+y_8 w_y))} \sinh \left(\frac{x}{4}(2\gamma+x_8 w_x+y_8 w_y)\right)	\times \\ 
		&\prod_{i=1}^{7} \sinh \left(\frac{x}{4}(2\gamma+x_i w_x+y_i w_y)\right)
	\end{align}

	Finally, specializing this expression to $sp$, we obtain 
	\begin{align}\label{uspq9-2} \nonumber
		&\frac{\sinh\left( \frac{x}{4}(2\gamma+\bar{x} \frac{3}{2}-\bar{y}\frac{1}{2})\right)}{ \sinh\left( \frac{x}{4}(2\gamma+\bar{x} \frac{3}{2}-y_8\frac{1}{2})\right) } \sinh \left(\frac{x}{4}(2\gamma+x_8 \frac{3}{2}-y_8 \frac{1}{2})\right) \times  \\
		&\prod_{i=1}^{7} \sinh \left(\frac{x}{4}(2\gamma+x_i \frac{3}{2}-y_i \frac{1}{2})\right) \\ \nonumber
		&\gamma=n+2
	\end{align}

	It is easy to see that the denominator can cancel neither the first hyperbolic sine in the numerator nor the second one (the latter involving $x_8$), since in either case the cancellation would occur identically. We therefore assume that it cancels, say, the seventh factor in the product, i.e. the one involving $x_7$. This means that 
	\begin{eqnarray}
		\bar{x} \frac{3}{2}-y_8 \frac{1}{2}	=x_7 \frac{3}{2}-y_7 \frac{1}{2}
	\end{eqnarray}
	
	The arguments of the remaining hyperbolic sines must then coincide with the same set of values as in (\ref{xyq8}):
	\begin{align}\label{xyq9} \nonumber
	&	\{	x_i \frac{3}{2}-y_i \frac{1}{2}, i=1...6, 	x_8 \frac{3}{2}-y_8 \frac{1}{2},	\bar{x} \frac{3}{2}-\bar{y} \frac{1}{2}	\}=   \\  
	&	\{-6,-4,-3,-2,-2,0,2,6\} \\
	&	\bar{x} \frac{3}{2}-y_8 \frac{1}{2}	=x_7 \frac{3}{2}-y_7 \frac{1}{2}
	\end{align}
	
	It is now easy to see that none of the $y_i$ for $i=1,...,6$ can be equal to 4, since there is no corresponding value of $x_i$ in (\ref{xq9}) such that $x_i \frac{3}{2}-y_i \frac{1}{2}$ belongs to the set in (\ref{xyq9}). Since there are two occurrences of 4 in (\ref{yq9}), we conclude that the only possible values of $x_i, y_i$ for $i=1,...,6$ are precisely those obtained above:
	\begin{eqnarray}\label{2xq8} \nonumber
	&	x_1=-8,x_2=8,x_3=-4,x_4=4, x_5=-2,x_6=2 \\
	&	y_1=-12, y_2=12, y_3=-4, y_4=8, y_5=0, y_6=6
	\end{eqnarray}
	
	We are then left with the system of equations
	\begin{align} \nonumber
		&y_7=\bar{y}=4, \{x_7,x_8\}=\{-2,2\}, \\ \nonumber
		&	\bar{x} \frac{3}{2}-y_8 \frac{1}{2}	=x_7 \frac{3}{2}-y_7 \frac{1}{2}, \\
		&x_8 \frac{3}{2}-y_8 \frac{1}{2}=-2,	\bar{x} \frac{3}{2}-\bar{y} \frac{1}{2}	=-2
	\end{align}
	From this system we obtain two solutions:
	\begin{eqnarray}
		\bar{x}=0, x_7=2, x_8=-2, y_8=-2
	\end{eqnarray}
	and
	\begin{eqnarray}
		\bar{x}=0, x_7=-2, x_8=2, y_8=10
	\end{eqnarray}
	
	For the first solution, the corresponding $2\gamma$-terms $Q_2$ are
	\begin{align}\label{u2g1} \nonumber
		Q_2&=\frac{	\sinh\left(\frac{x}{4}(2\gamma+2\beta+2\alpha)\right)\sinh\left(\frac{x}{4}(2\gamma+\beta)\right) \sinh\left(\frac{x}{4}(2\gamma+\alpha)\right) }{ \sinh\left(\frac{x}{4}(2\gamma-\beta-\alpha)\right) } \times \\ \nonumber
	&	\sinh\left(\frac{x}{4}(2\gamma+2\beta-2\alpha)\right)\sinh\left(\frac{x}{4}(2\gamma-2\beta+2\alpha)\right) \times \\ \nonumber
	&	\sinh\left(\frac{x}{4}(2\gamma+2\beta+\alpha)\right)\sinh\left(\frac{x}{4}(2\gamma+\beta+2\alpha)\right) \times \\
	&	\sinh\left(\frac{x}{4}(2\gamma+2\beta)\right)\sinh\left(\frac{x}{4}(2\gamma+2\alpha)\right)
	\end{align}

	For the second solution, the $2\gamma$-terms $Q_2'$ are

	\begin{align}\label{u2g2} \nonumber
		Q_2' & =	\frac{	\sinh\left(\frac{x}{4}(2\gamma+2\beta+2\alpha)\right)\sinh\left(\frac{x}{4}(2\gamma+4\beta +3\alpha)\right) }{ \sinh\left(\frac{x}{4}(2\gamma+5\beta+5\alpha)\right) } \times \\ \nonumber
	&	\sinh\left(\frac{x}{4}(2\gamma+3\beta+4\alpha)\right) 	\sinh\left(\frac{x}{4}(2\gamma+2\beta-2\alpha)\right)\times \\ \nonumber
	&\sinh\left(\frac{x}{4}(2\gamma-2\beta+2\alpha)\right) 	\sinh\left(\frac{x}{4}(2\gamma+2\beta+\alpha)\right) \times \\ 
	&	\sinh\left(\frac{x}{4}(2\gamma+\beta+2\alpha)\right) \sinh\left(\frac{x}{4}(2\gamma+2\beta)\right)\sinh\left(\frac{x}{4}(2\gamma+2\alpha)\right)
	\end{align}

	It remains to determine the $1\gamma$-terms in the universal quantum dimension formula for the universal multiplet $E$.

	\section{Guess of the final formula}
	
	At this stage, we can already guess the final formula. So far, we have obtained the $\gamma$-independent factor (\ref{udenomE2}) and the $2\gamma$-terms, given by either (\ref{u2g1}) or (\ref{u2g2}). The final formula should be the product of these factors together with the as yet undetermined $1\gamma$-terms. Consider the product of (\ref{udenomE2}) and the first solution (\ref{u2g1}) in the limit $x\rightarrow 0$, and compare it with the dimension formula (\ref{dimE}). One then has to include the factors from the dimension formula that are missing in this product, and at the same time introduce factors cancelling those that are present in the product but absent from the dimension formula. This leads to the following $1\gamma$-terms in the quantum dimension:

	\begin{align}\label{1gterm} \nonumber
		Q_1=	&\frac{	\sinh\left(\frac{x}{4}(\gamma+2\alpha)\right) 	\sinh\left(\frac{x}{4}(\gamma+2\beta)\right) 	\sinh\left(\frac{x}{4}(\gamma+2\beta+\alpha)\right) 	}{	\sinh\left(\frac{x}{4}\gamma\right)	\sinh\left(\frac{x}{4}(\gamma-\alpha)\right)	\sinh\left(\frac{x}{4}(\gamma-\beta)\right) 	 }
		\times \\ 
		& \frac{		\sinh\left(\frac{x}{4}(\gamma+\beta+2\alpha)\right)	\sinh\left(\frac{x}{4}(\gamma+2\beta+2\alpha)\right) }{		\sinh\left(\frac{x}{4}(\gamma+\beta-\alpha)\right) 	\sinh\left(\frac{x}{4}(\gamma-\beta+\alpha)\right) }	
	\end{align}

	Altogether, the final universal quantum dimension formula for the $E$ multiplet is the product of (\ref{udenomE2}), (\ref{u2g1}), and (\ref{1gterm}):
	
	\begin{eqnarray}\label{FinalE}
		\text{Universal quantum dimension of } E=Q_0Q_1Q_2
	\end{eqnarray}
	
	This formula can be verified for all algebras, as well as for permutations of the parameters, and in each case it gives the correct quantum dimensions. In the next section we present explicit formulae for comparison. The second solution (\ref{u2g2}) does not yield the correct answers.
	
	\section{Quantum dimensions of $E$ for exceptional algebras}
	
	For the five exceptional algebras, the formula (\ref{FinalE}) gives the following quantum dimensions, which agree with direct calculations in the corresponding algebras using the Weyl character formula. 
	
	For the algebra $E_8$, the dimension formula (\ref{dimE}) gives 26411008, i.e. the dimension of the representation with Dynkin labels (01000001) in Bourbaki's \cite{Bourbaki1981} numbering (which will be used for all exceptional algebras everywhere below; we also do not put comma between different labels, since in this paper they are always one-digit). Its quantum dimension is given below and coincides with (\ref{FinalE}) specialized to the $E_8$ parameters from table \ref{tab:Vogel}:

	\begin{align}\label{E8} \nonumber
		&\frac{
			\sinh\!\left(\frac{8x}{2}\right)
			\sinh\!\left(\frac{14x}{2}\right)
			\sinh\!\left(\frac{18x}{2}\right)
			\sinh\!\left(\frac{20x}{2}\right)
			\sinh\!\left(\frac{21x}{2}\right)
			\sinh\!\left(\frac{22x}{2}\right)
		}{
			\sinh^2\!\left(\frac{x}{2}\right)
			\sinh\!\left(\frac{3x}{2}\right)
			\sinh(2x)
			\sinh\!\left(\frac{5x}{2}\right)
			\sinh(3x)
			\sinh^2\!\left(\frac{7x}{2}\right)
		}  
		\times \\ 
	&	\frac{ 	\sinh\!\left(\frac{24x}{2}\right)
			\sinh\!\left(\frac{26x}{2}\right)
			\sinh\!\left(\frac{30x}{2}\right)
			\sinh\!\left(\frac{31x}{2}\right)
			\sinh\!\left(\frac{32x}{2}\right)
			\sinh\!\left(\frac{34x}{2}\right)}{ 	
			\sinh\!\left(\frac{9x}{2}\right)
			\sinh\!\left(\frac{11x}{2}\right)
			\sinh\!\left(\frac{15x}{2}\right)
			\sinh\!\left(\frac{17x}{2}\right)
		}
	\end{align}

	For $E_7$ we obtain the representation (1001001), with dimension 3424256 and quantum dimension

	\begin{align}\label{E7} \nonumber
		&\frac{
			\sinh\!\big(\tfrac{8x}{2}\big)\,
			\sinh\!\big(\tfrac{10x}{2}\big)\,
			\sinh\!\big(\tfrac{12x}{2}\big)\,
			\sinh\!\big(\tfrac{14x}{2}\big)^{2}\,
			\sinh\!\big(\tfrac{16x}{2}\big)\,	
		}{
			\sinh\!\big(\tfrac{x}{2}\big)^{3}\,
			\sinh\!\big(\tfrac{3x}{2}\big)\,
			\sinh\!\big(\tfrac{4x}{2}\big)\,
		}   
		\times \\	
		&\frac{ 	\sinh\!\big(\tfrac{18x}{2}\big)\,
			\sinh\!\big(\tfrac{19x}{2}\big)\,
			\sinh\!\big(\tfrac{20x}{2}\big)\,
			\sinh\!\big(\tfrac{22x}{2}\big) }{ 
			\sinh\!\big(\tfrac{5x}{2}\big)^{2}\,
			\sinh\!\big(\tfrac{7x}{2}\big)^{2}\,
			\sinh\!\big(\tfrac{9x}{2}\big) }
	\end{align}

	For $E_6$, formula (\ref{dimE}) gives the dimension 504504, which is twice the dimension of the representation (111000). This is due to the $Z_2$ symmetry of the Dynkin diagram. The corresponding quantum dimension is twice the following expression:

	\begin{align}\label{E6} \nonumber
	&	\frac{
			\sinh\!\big(\frac{7x}{2}\big)\,
			\sinh\!\big(\frac{9x}{2}\big)^{2}\,
			\sinh\!\big(\frac{11x}{2}\big)\,
			\sinh\!\big(\frac{12x}{2}\big)
		}{
			\sinh\!\big(\frac{x}{2}\big)^{3}\,
			\sinh\!\big(\frac{3x}{2}\big)^{2}\,
			\sinh\!\big(\frac{4x}{2}\big)^{2}\,
			\sinh\!\big(\frac{6x}{2}\big)}  \times  \\			
			&	\sinh\!\big(\frac{13x}{2}\big)\,
				\sinh\!\big(\frac{14x}{2}\big)\,
				\sinh\!\big(\frac{16x}{2}\big)
	\end{align}

	The universal formula (\ref{FinalE}) produces the coefficient 2 due to the factor (the hyperbolic sine in the numerator comes from $Q_2$, while that in the denominator comes from $Q_1$):
	
	\begin{eqnarray}
		\frac{\sinh\left(\frac{x}{4}(2\gamma-2\beta+2\alpha)\right)}{\sinh\left(\frac{x}{4}(\gamma-\beta+\alpha)\right)}
	\end{eqnarray}
	which, for the $E_6$ parameters $\alpha=-2, \beta=6, \gamma=8$, tends to 2, while the remaining part of $Q_0Q_1Q_2$ gives (\ref{E6}).
	
	For $F_4$ we obtain, in both ways, the representation (1011), with dimension 106496 and quantum dimension:
	
	\begin{align}\label{F4} \nonumber
		&\frac{
			\sinh\!\left(\frac{5x}{2}\right)\,
			\sinh(3x)\,
			\sinh\!\left(\frac{7x}{2}\right)\,
			\sinh(4x)\,		
		}{
			\sinh^2\!\left(\frac{x}{4}\right)\,
			\sinh\!\left(\frac{3x}{4}\right)\,
			\sinh^2\!\left(\frac{5x}{4}\right)\,		
		} \times \\
		&\frac{	\sinh\!\left(\frac{9x}{2}\right)\,
			\sinh(5x)\,
			\sinh\!\left(\frac{11x}{2}\right)\,
			\sinh\!\left(\frac{13x}{2}\right)
		}{
			\sinh\!\left(\frac{7x}{4}\right)\,
		\sinh\!\left(\frac{9x}{4}\right)\,
		\sinh\!\left(\frac{11x}{4}\right)
		}
	\end{align}

	For $G_2$ we obtain the representation (12), with dimension 286 and quantum dimension:
	
	\begin{eqnarray}\label{G2}
		\frac{
			\sinh\!\left(\frac{x}{3}\right)
			\sinh\!\left(\frac{11x}{6}\right)
			\sinh\!\left(\frac{13x}{6}\right)
			\sinh\!\left(\frac{5x}{2}\right)
			\sinh(4x)
		}{
			\sinh\!\left(\frac{x}{6}\right)
			\sinh\!\left(\frac{x}{2}\right)
			\sinh\!\left(\frac{2x}{3}\right)
			\sinh\!\left(\frac{5x}{6}\right)
			\sinh(x)
		}
	\end{eqnarray} 
	and the same answer from universal formula (\ref{FinalE}).
	
	We also checked the main formula for the permutations of parameters for the $E_8$ algebra. There are three such permutations, since the formula is symmetric under $\alpha \leftrightarrow \beta$. For $\alpha=12, \beta=20, \gamma=-2$, the universal dimension formula gives a virtual (i.e. taken with a minus sign) representation with Dynkin labels (00100000), dimension 147250, and quantum dimension
	\begin{align}\label{E8b} \nonumber
		&\frac{
			\sinh\!\left(\frac{19x}{2}\right)
			\sinh(10x)
			\sinh\!\left(\frac{21x}{2}\right)
			\sinh\!\left(\frac{25x}{2}\right)		
		}{
			\sinh\!\left(\frac{x}{2}\right)
			\sinh\!\left(\frac{3x}{2}\right)
			\sinh(2x)
			\sinh\!\left(\frac{5x}{2}\right)		
		} \times \\
		&\frac{	\sinh(15x)
			\sinh\!\left(\frac{31x}{2}\right)
			\sinh(16x)
		}{
			\sinh(3x)
		\sinh\!\left(\frac{7x}{2}\right)
		\sinh(8x)
		}
	\end{align}

	The universal formula (\ref{FinalE}), in complete agreement, yields the same expression with a minus sign. 
	
	Similarly, for $\alpha=-2, \beta=20, \gamma=12$, the universal dimension formula gives a virtual (i.e. again with a minus sign) representation with Dynkin labels (00000011), dimension 4096000, and quantum dimension
	\begin{align}\label{E8g} \nonumber
		&\frac{
			\sinh\!\left(\frac{14x}{2}\right)
			\sinh\!\left(\frac{18x}{2}\right)
			\sinh\!\left(\frac{20x}{2}\right)
			\sinh\!\left(\frac{22x}{2}\right)
			\sinh\!\left(\frac{24x}{2}\right)
		}{
			\sinh^2\!\left(\frac{x}{2}\right)
			\sinh\!\left(\frac{3x}{2}\right)
			\sinh\!\left(\frac{5x}{2}\right)
			\sinh\!\left(\frac{6x}{2}\right)
			\sinh\!\left(\frac{7x}{2}\right)	
		} \times \\ 
	&	\frac{  \sinh\!\left(\frac{25x}{2}\right)
			\sinh\!\left(\frac{26x}{2}\right)
			\sinh\!\left(\frac{30x}{2}\right)
			\sinh\!\left(\frac{32x}{2}\right)
			\sinh\!\left(\frac{34x}{2}\right)  }{  
			\sinh\!\left(\frac{9x}{2}\right)
			\sinh\!\left(\frac{11x}{2}\right)
			\sinh\!\left(\frac{13x}{2}\right)
			\sinh\!\left(\frac{17x}{2}\right) }
	\end{align}

	Again, the universal formula (\ref{FinalE}), in complete agreement, gives the same expression with a minus sign.

	\section{$1\gamma$-terms from quantum dimensions of algebras}
	
	How can $1\gamma$ terms be obtained in a systematic way? Examining the quantum dimensions of the $sl$, $so$, and $sp$ algebras, we see that the minimal number of such terms is two in the numerator and two in the denominator. One can easily find a possible expression that gives the required values for all three algebras:
	\begin{align} \nonumber
	&	\frac{\sinh\left(\frac{x}{4}(\gamma+4w_y)\right) \sinh\left(\frac{x}{4}(\gamma+8w_y)\right)}{\sinh\left(\frac{x}{4}(\gamma-6w_y)\right) \sinh\left(\frac{x}{4}(\gamma+2w_y)\right)} = \\ 
	&	\frac{\sinh\left(\frac{x}{4}(\gamma+2\beta+2\alpha)\right) \sinh\left(\frac{x}{4}(\gamma+4\beta+4\alpha)\right)}{\sinh\left(\frac{x}{4}(\gamma-3\beta-3\alpha)\right) \sinh\left(\frac{x}{4}(\gamma+\beta+\alpha)\right)}
	\end{align}
	
	Indeed, for $sl$ this expression is equal to 1, while for $so(N)$ it becomes 
	\begin{eqnarray}
		\frac{\sinh\left(\frac{x}{4}N \right)
			\sinh\left(\frac{x}{4}(N+4)\right)}{\sinh\left(\frac{x}{4}(N-10)\right) \sinh\left(\frac{x}{4}(N-2)\right)}
	\end{eqnarray}
	in agreement with (\ref{qsoe}), and for $sp(2n)$ it becomes 
	\begin{eqnarray}
		\frac{\sinh\left(\frac{x}{4}n \right)
			\sinh\left(\frac{x}{4}(n-2)\right)}{\sinh\left(\frac{x}{4}(n+5)\right) \sinh\left(\frac{x}{4}(n+1)\right)}
	\end{eqnarray}
	in agreement with (\ref{qspe}). However, this universal expression does not reproduce the quantum dimensions for the exceptional algebras; in fact, it already fails to give the correct universal dimension formula (\ref{dimE}).
	
	We therefore seek to use additional information, namely that contained in the expressions for the quantum dimensions of the exceptional algebras, in order to reconstruct the universal expression for the quantum dimensions of $E$. 
	
	For $E_8$, one can determine the value of the $1\gamma$-terms simply by dividing its quantum dimension by the product of the $\gamma$-independent factor and the $2\gamma$-terms. 
	
	Before doing so, however, we must determine which representations of the exceptional algebras enter the universal multiplet $E$. The point is that, although the representation can be identified from the universal dimension formula (\ref{dimE}), our aim is to derive new universal dimension formulae, and therefore in general we do not know in advance which representations belong to the corresponding multiplet. 
	
	Above, we started from a representation of the $sl$ algebra corresponding to two Young diagrams $\lambda$ and $\mu$ of equal area, and constructed, according to \cite{Mkrtchyan2025b,Mkrtchyan2026}, the associated representations of $so$ and $sp$ as those with Young diagrams given by the vertical sum $\lambda \oplus_v \mu$ and the horizontal sum $\lambda \oplus_h \mu$, respectively. We now need to extend this correspondence to exceptional algebras. 
	
	Although a general recipe is not yet known, in the present case, as well as in a number of similar cases, this can be done in the following way. The nonzero representations of the universal multiplet $E$ in the small universal multiplets of classical algebras are given by the Cartan product of the adjoint representation with another representation, namely the permuted Cartan product of the adjoint representation with $X_2$. In Vogel's notation this representation is denoted by $C'$, while in the notation of \cite{LandsbergManivel2006} it is denoted by $C(\beta)$. These representations admit a universal description, and their analogues for the exceptional case can therefore be identified. In particular, for $so$ the representation $C(\beta)$ has Dynkin labels $(100010...)$, and its Cartan product with the adjoint representation has labels $(110010...)$, in agreement with (\ref{soE}). In the case of the $E_8$ algebra, $C(\beta)$ is the fundamental representation corresponding to root number 2 in Bourbaki's numbering. Since the adjoint representation is the fundamental representation corresponding to root number 8, the desired representation of $E_8$ has highest weight $w_2+w_8$, again in agreement with (\ref{E8}). For the other exceptional algebras, the highest weights of E can be found in a similar way and are listed in the previous section. 
	
	It is also worth mentioning that we know the universal representations corresponding to at least two other fundamental weights of $E_8$: the weight $w_7$ corresponds to the representation $X_2$, and the weight $w_1$ corresponds to the representation $Y_2(\beta)$. Since the number of fundamental weights of $E_8$ is eight, and we already know the universal representations corresponding to four of them, it seems plausible that one can establish a universal interpretation for all representations of $E_8$, and in particular determine their correspondence with the associated representations of classical algebras. For the other exceptional algebras this problem should be easier, since their rank, i.e. the number of fundamental weights, is smaller. On the other hand, we expect that some representations of these algebras (for example, defining representations or those of smallest dimension) cannot be characterized universally, as is the case for the defining representations of the classical algebras. Note that $E_8$ is the unique algebra whose smallest nontrivial representation coincides with the adjoint representation. We shall return to these questions elsewhere. 
	
	We therefore determine the $1\gamma$-terms for $E_8$ as described above:	
	\begin{eqnarray}\label{1gterms}
		\frac{\sinh\left(\frac{x}{4}16\right) \sinh\left(\frac{x}{4}28\right) \sinh\left(\frac{x}{4}40\right) \sinh\left(\frac{x}{4}42\right)\sinh\left(\frac{x}{4}44\right)}{   \sinh\left(\frac{x}{4}6\right)  \sinh\left(\frac{x}{4}8\right) \sinh\left(\frac{x}{4}20\right) \sinh\left(\frac{x}{4}22\right)\sinh\left(\frac{x}{4}34\right) }
	\end{eqnarray}
	
	We see that this is a ratio of five hyperbolic sines to five hyperbolic sines, exactly as in the universal $1\gamma$-term (\ref{1gterm}) obtained earlier. Our aim is to rederive independently the universal form (\ref{1gterm}) of the $1\gamma$-terms from (\ref{1gterms}) and the remaining data. 
	
	From the $E_8$ expression (\ref{1gterms}) we see that the desired universal form should involve at least five hyperbolic sines in the numerator and five in the denominator, with universal arguments. We assume this minimal number and consider its most general form:
	\begin{eqnarray}\label{u5over5}
		\prod_{i=1}^{5} \frac{\sinh\left(\frac{x}{4}(\gamma+x_iw_x+y_iw_y)\right) }{ \sinh\left(\frac{x}{4}(\gamma+\bar{x}_i w_x+\bar{y}_i w_y)\right) }
	\end{eqnarray}
	
	For $sl$ this expression must be equal to 1, since there are no such terms in (\ref{qsle}). This implies that, after a suitable reindexing, $\bar{x}_i=x_i$ for $i=1,...,5$. When specialized to the $E_8$ case, the arguments of the five hyperbolic sines in the numerator must coincide with those in (\ref{1gterms}), and similarly for the denominator. This observation allows one to express ten of the fifteen variables $x_i$, $y_i$, and $\bar{y}_i$ in terms of the remaining five, provided that the arguments are ordered increasingly. 
	
	For $so$, three pairs of hyperbolic sines must cancel, which means that three of the $\bar{y}_i$ must coincide with three of the $y_i$. One can list all possible cancellation patterns for the $so$ algebra. For each such choice one obtains an overdetermined system of linear equations. These systems can be solved straightforwardly because they are linear, although a solution need not exist in general since the system is overdetermined. In fact, the values of the arguments of only two hyperbolic sines in the numerator for $so$ already fix all variables, so that the arguments of the hyperbolic sines in the denominator for $so$, as well as the full set of $1\gamma$-terms for $sp$, should then emerge automatically. We omit the rather lengthy case-by-case analysis and instead anticipate the simple proof of uniqueness in terms of the theory of configurations of points and lines \cite{AvetisyanMkrtchyan2022}. The actual solution arises when one chooses the cancellation pattern $y_1=\bar{y}_2$, $y_2=\bar{y}_3$, and $y_4=\bar{y}_5$. Then the parameters in (\ref{u5over5}) are given by 
	\begin{align} \nonumber
		&x_1=-4,\quad x_2=-2, \quad x_3=0, \quad x_4=2, \quad x_5=4, \\ \nonumber
		&y_1=-4,\quad y_2=0, \quad y_3=4, \quad y_4=6, \quad y_5=8, \\
		&\bar{y}_1=-6,\quad \bar{y}_2=-4, \quad \bar{y}_3=0, \quad \bar{y}_4=2, \quad \bar{y}_5=6 
	\end{align}
	which, after substitution into (\ref{u5over5}), yields (\ref{FinalE}), thus completing the derivation.

	\section{Conclusion}
	
	In this paper, together with the first part \cite{Mkrtchyan2026}, we have proposed a method for calculating universal quantum dimensions for a certain class of universal multiplets. We applied this method to the universal multiplet E and obtained a new formula for its universal quantum dimension.
	
	Let us now summarize, in general form, all steps involved in the calculation, emphasizing the ambiguities and choices that arise in the process. 
	
	\begin{enumerate}
		\item Choose a representation of $sl$ that will serve as a member of the universal multiplet. This representation should belong to the set of representations appearing in the decomposition of $\mathfrak{g}^k$, with the areas of its constituent $\lambda$ and $\mu$ representations equal to $k$. It should also appear ``alone'' in the sense that there is no other representation in the decomposition of $\mathfrak{g}^k$ with the same Casimir eigenvalue (apart from the symmetrized one obtained by interchanging $\lambda$ and $\mu$). Otherwise, one may be dealing not with a universal multiplet, but with a universal Casimir multiplet, i.e. the sum of irreducible representations with the same Casimir eigenvalue, with universal expression for the sum of their dimensions.
		
		\item Determine the associated $so$ and $sp$ members, which can be obtained straightforwardly as the vertical and horizontal sums of the Young diagrams corresponding to $\lambda$ and $\mu$. 
		
		\item  One may need the associated representations of the exceptional algebras. The rules for finding these representations are currently known only partially. In particular, if, say, the $so$ member of the multiplet is given by the Cartan product of representations that already admit a universal description, then the corresponding representation of an exceptional algebra should be given by the same Cartan product of the corresponding representations of exceptional. Thus, it remains to determine which universal representations correspond to the (finite) set of fundamental representations of the exceptional algebras. 
		
		\item Compute the quantum dimensions of all these representations. This is a standard Lie-algebraic calculation. 
		
		\item Choose the pattern, i.e. the number of hyperbolic sines in the numerator and denominator, for the $\gamma$-independent, $2\gamma$-, and $1\gamma$-terms. This is a discrete choice, constrained from below by the explicit expressions for the quantum dimensions of the algebras. One starts with the minimal possibility and, if it does not work, increases the numbers of hyperbolic sines simultaneously in the numerator and denominator. Note the subtleties coming from the automorphisms of Dynkin diagrams: they  lead  to seemingly singular terms, which actually have a finite limit. 
		
		\item For a given pattern  compare expression obtained with the results for the actual algebras and determine the universal formulae for $\gamma$-independent, $1\gamma$- and $2\gamma$-terms (i.e. the numbers $x_i, y_i, \bar{x}_i, \bar{y}_i$ in above). We assume that either this will prove impossible, in which case one has to pass to the next pattern, or else it will yield the correct answer. We regard this step as the core of the method. It replaces complicated calculations of other methods with completely different calculations and seems much more affordable.
		
		\item The final result must be checked completely, particularly for the permutations of parameters. 
		
	\end{enumerate}
	
	Some of these steps are straightforward, others are more complicated, and some of them may require repetition with new sets of parameters. We expect that this procedure will make it possible to find new universal multiplets, if not all in the scope, but at least those appearing in sufficiently low powers (fifth and higher) of the adjoint representation. Moreover, the algorithm is well suited for computer implementation.
	
	One can consider the extension of this method to the refined version of the universality, see \cite{KreflSchwarz2013,AvetisyanMkrtchyan2021,BishlerMironovMorozov2025,BishlerMironov2025b,Bishler2025,MironovSingh2025}. All steps above can be carried on for refined case, too. The possible problem may be the different behavior under universalization of simply-laced and non-simply-laced algebras, so perhaps first simply-laced algebras should be considered.

	\section*{Acknowledgments}
	
	This work was partially supported by the Science Committee of the Ministry of Science and Education of the Republic of Armenia under contracts 21AG-1C060 and 24WS-1C031.

	\appendix
	
	\section{Vogel's table}
	
	Here we reproduce  Vogel's table \ref{tab:Vogel} of points in the projective plane corresponding to simple Lie algebras. 
	\begin{table}[h] \caption{Vogel's parameters for simple Lie algebras}     \label{tab:Vogel}
		\begin{tabular}{|r|r|r|r|r|r|} 
			\hline & $\alpha$ &$\beta$  &$\gamma$  & $t=  \alpha+\beta+\gamma$ & Line \\ 
			\hline $sl(N)$ & -2 & 2 & $N$ & $N$ & $\alpha+\beta=0 $\\ 
			\hline $so(N) $ & -2  & 4 & $N-4$ & $N-2$ & $2\alpha+\beta=0$ \\ 
			\hline $sp(2n)$ & -2  & 1 & $n+2$ & $n+1$ & $\alpha+2\beta=0$ \\ 
			\hline $Exc(k)$ & -2 & $k+4$  & $2k+4$ & $3k+6$& $\gamma=2(\alpha+\beta)$ \\ 
			\hline 
		\end{tabular} 
	\end{table}
	
	In the table \ref{tab:Vogel} for the exceptional line $Exc(k)$, $k=-1,-2/3,0,1,2,4,8$ correspond to $A_2,G_2, D_4, F_4, E_6, E_7, E_8$, respectively.

	
	
	

	
\end{document}